\documentclass[12pt, reqno]{amsart}
\usepackage{amssymb, amsthm, amsmath, amsfonts}
\usepackage{array, texdraw}
\evensidemargin \oddsidemargin
\begin{document}

\renewcommand{\proofname}{\bf Proof}
\newtheorem*{rem*}{Remark}
\newtheorem*{cor*}{Corollary}
\newtheorem{pred}{Proposition}
\newtheorem{lem}{Lemma}
\newtheorem{theo}{Theorem}

\def\R{\mathbb{R}}
\def\N{\mathbb{N}}
\def\E{\mathbb{E}}
\def\P{\mathbb{P}}
\def\D{\mathbb{D}}
\def\I{\mathbf{1}}

\title[A limit theorem for the position of a
particle in a Lorentz type model]{A functional limit theorem for
the\\ position of a particle in a Lorentz type model}

\author[V. Vysotsky]{Vladislav V. Vysotsky}

\begin{abstract}
Consider a particle moving through a random medium, which consists
of spherical obstacles, randomly distributed in $\R^d$. The
particle is accelerated by a constant external field; when
colliding with an obstacle, the particle inelastically reflects.
We study the asymptotics of $X(t)$, which denotes the position of
the particle at time $t$, as $t \to \infty$. The result is a
functional limit theorem for $X(t)$.

{\it Key words and phrases:} Lorentz model, motion in random
medium, functional central limit theorem for Markov chains, limit
theorems.

{\it 2000 Mathematics Subject Classification:} {\bf 60K37}
\end{abstract}


\maketitle

\section{The Lorentz model and the problem}
\subsection{The motivation} \label{Motiv}
Consider a spherical particle moving in a random medium. The
medium consists of immobile spherical obstacles of equal radii,
randomly distributed in $\R^3$. The particle is accelerated by an
external field providing constant acceleration $a$. At a collision
with an obstacle, the particle's speed $v$ changes to
$$v - (1+\alpha)(v, \nu) \nu,$$ where $0 \le \alpha \le 1$ is the
restitution coefficient and $\nu$ is the inner unit normal to the
obstacle at the point of collision.
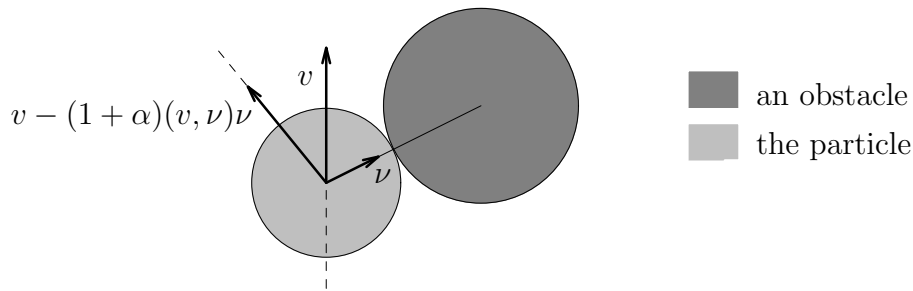
\begin{center}
\begin{texdraw}
\drawdim cm \move(1.2 1.2) \linewd 0.01 \fcir f:0.75 r:1
\move(3.257 2.229) \fcir f:0.5 r:1.3 \linewd 0.01 \move(1.2 1.2)
\lcir r:0.99 \move(3.257 2.229) \lcir r:1.299 \move(1.2 1.2)
\move(1.2 1.2) \linewd 0.003 \lvec(3.257 2.229) \linewd 0.012

\lpatt(0.12 0.12) \move(1.2 -0.2) \lvec(1.2 1.2) \lvec(-.225
2.954) \lpatt ( )

\move(1.2 1.2) \linewd 0.035 \arrowheadsize l:0.2 w:0.1
\arrowheadtype t:V \avec (1.2 3) \move(1.2 1.2) \avec(1.9 1.55)
\move(1.2 1.2) \avec (0.16 2.48)

\move(1.02 2.6) \textref h:R v:C \htext{$v$} \move(2.05 1.3)
\textref h:R v:C \htext{$\nu$} \move(0.2 2.1) \textref h:R v:C
\htext{$v - (1+\alpha)(v, \nu) \nu$}

\setgray 1 \move(6 2.2) \lvec(6 2.7) \lvec(6.7 2.7) \lvec(6.7 2.2)
\lvec(6.7 2.2) \lfill f:0.5 \move(6.9 2.4) \textref h:L v:C
\htext{an obstacle}

\move(6 1.5) \lvec(6 2) \lvec(6.7 2) \lvec(6.7 1.5) \lvec(6.5 1.5)
\lfill f:0.75 \move(6.9 1.7) \textref h:L v:C \htext{the particle}

\end{texdraw}
\end{center}
This mapping changes only the normal component of $v$, i.e., $(v,
\nu) \nu$, which is multiplied by $-\alpha$. On the figure, the
dotted line indicates the trajectory of the particle's center.

This model is often named after Hendrik Lorentz who introduced it
(see \cite{Lor}) in order to describe conductivity in metals.
Lorentz studied the case of elastic collisions, with $\alpha=1$;
the generalization for non-elastic collisions could be found,
e.g., in \cite{WE}.

In physics the Lorentz model is used to describe the motion of a
particle in a medium if the particle's mass is negligible with
respect to masses of the medium's particles (obstacles). Indeed,
in this case we can assume that the obstacles have infinite masses
and thus they remain immobile at collisions. For example, the
model in question describes well the motion of electrons in
helium, see \cite{Helium}.

In case of elastic collisions, the Lorentz model is a billiard
type model; recall that {\it a billiard} is the following model of
motion of a particle in an arbitrary region with smooth boundary:
at a collision with the boundary, the particle elastically
reflects. Indeed, we can consider the region whose boundary
consists of the obstacles' boundaries. The case when the obstacles
are located periodically, is especially interesting (see
\cite{BSinai}, \cite{San}); such billiards are called {\it Lorentz
periodic gases}. The main tool for studying billiards is ergodic
theory; the basic results of this theory and their applications
for the Lorentz model could be found in \cite{Erg}. For a detailed
review of problems and methods of the billiard theory, see
\cite{Erg} and \cite{Hard}.

Another interesting interpretation of the Lorentz model could be
found in \cite{WE} (also, see the references therein): a particle
{\it percolates} through an immobile medium under the constant
gravity. Considering a big number of such percolating particles
and neglecting interactions between them, we get a model of mixing
of two dry substances, for example, powders.

There are many physical papers on the Lorentz model. Usually their
main goal is to study the Boltzmann equation (that is the equation
for density $p(v, t)$ of probability that at time $t$ the
particle's speed is $v$; this equation could be derived from the
law of conservation of matter); see \cite{BFS}, \cite{BCLM}, and
\cite{MP}. For instance, the purpose of \cite{MP} is to prove the
existence of the stationary, i.e., independent of $t$, solution of
the Boltzmann equation in case $\alpha < 1$. It is typical that in
\cite{MP} the convergence of $p(v, t)$ to the stationary solution,
as $t \to \infty$, is not discussed.

\subsection{The model} \label{my model}
Since the Lorentz model is extremely complicated to analyze, we
replace it with a simpler one. In doing so, we follow multiple
papers on the Lorentz model, for example, \cite{BFS}, \cite{BCLM},
\cite{MP}, and \cite{WE}. For the case of elastic collisions, our
simplified model coincide with one introduced in \cite{R-Triolo};
in this paper the authors only consider $\alpha = 1$ and formulate
their model in a rather different way. We also note that the
present simplified model was implicitly used in \cite{MP}.
However, in \cite{MP} and \cite{R-Triolo} the authors do not
discuss how the simpler model is derived from the original one. In
Section~\ref{Append} we give such explanation.

Let us formulate the simplified model. Denote by $V_n$ the speed
of the particle just before the $n$th collision and denote by
$\tau_n$ the random time between the $n$th and the $(n+1)$th
collisions. Let $\{\sigma_n\}_{n \ge 1} \subset S^2 \subset \R^3$
be uniformly distributed unit vectors; let $\{ \eta_n \}_{n \ge
0}$ be exponential random variables with mean $\lambda$; and let
$\{\sigma_n\}_{n \ge 1}, \{ \eta_n \}_{n \ge 0}$ be independent.
Here $\lambda > 0$ is a parameter signifying the mean free path of
the particle. We state that
\begin{equation}
\label{V_n+1 long} V_{n+1} = V_n - \frac{1 + \alpha}{2} \bigl(V_n
+ |V_n| \sigma_n \bigr) + a \tau_n;
\end{equation}
\begin{equation}
\label{tau_n def} \tau_n=F \Bigl(V_n - \frac{1 + \alpha}{2}
\bigl(V_n + |V_n| \sigma_n \bigr), \eta_n \Bigr),
\end{equation}
where $F: \R^3 \times \R_+ \to \R$ is a deterministic function
defined as the solution of the equation $$\int_0^{F(v, t)} |v +
as|ds = t.$$ In addition, we suppose that at time zero the
particle's speed is nonrandom and equals some $v_0 \in \R^3$.
Therefore $V_1 = v_0 + a \tau_0,$ where $\tau_0=F (v_0, \eta_0)$
is the random moment of the first collision. We stress that $V_n$
form a Markov chain, and the motion of the particle is completely
defined by this chain.

Let us briefly compare the new model and the original one. In the
new model, at the $n$th collision the speed $V_n$ changes to $V_n
- \frac{1 + \alpha}{2} (V_n + |V_n| \sigma_n )$. This corresponds
to the collision with an obstacle, whose inner normal $\nu_n$
directed along the bisectrix of the angle between $\sigma_n$ and
$V_n$ (in fact, $(V_n, \nu_n) \nu_n = \frac{1}{2}(V_n + |V_n|
\sigma_n)$). Further, since for every $v$ the function $F(v,
\cdot)$ is monotone, we have $F \left(v, \int_0^t |v + as|ds
\right) = t$. Substituting $v:=V_n - \frac{1 + \alpha}{2} (V_n +
|V_n| \sigma_n )$ and $t:=\tau_n$, and comparing the resulting
equality with \eqref{tau_n def}, we get
$$\eta_n = \int_0^{\tau_n} \Bigl|V_n - \frac{1 + \alpha}{2}
(V_n + |V_n| \sigma_n ) + as \Bigr|ds.$$ Thus $\eta_n$ is the
length of the path passed by the particle between the $n$th and
the $(n+1)$th collisions. As the mean free path is $\lambda$, it
is quite natural that $\eta_n$ are exponential r. vs. with mean
$\lambda$. If we know the path's length $\eta_n$ and the initial
speed $V_n - \frac{1 + \alpha}{2} (V_n + |V_n| \sigma_n )$, we can
find the time $\tau_n$; this argument explains \eqref{tau_n def}.

However, there are quite significant disparities between the
models. In the new model, if the particle collides with an
obstacle at some point of space, then the particle does not
necessarily collide again while coming back to the same point. In
other words, after a collision the obstacle instantly
"disappears". This happens because of the postulated independence
of all $\sigma_n, \eta_n$. Thus the medium "changes" in a rather
specific way. See more about the disparities in
Section~\ref{Append}.

\subsection{The problem and the results}
We study the asymptotics of $X(t)$, which denotes the position of
the particle at time $t$, as $t \to \infty$. The case $\alpha = 1$
was investigated in \cite{R-Triolo}, where the authors proved that
after a proper normalization the trajectories of $X(t)$ weakly
converge to a certain diffusion process.

Let us mention about a very similar model of motion in $\R^1$: at
a collision, the particle's speed $v$ always changes to $-\alpha
v$. If there is no external field and $\alpha = 1$, then the
times between collisions are i.i.d. exponential r.vs., thus $X(t)$
is a well known telegraph process. This simplest model of motion
is studied in detail in \cite{Kac}.

In the current paper, we consider the motion in $\R^3$, with an
external field, and assume that $\alpha \in (0,1)$. The purpose of
this paper is to prove a functional limit theorem for $X(t)$, and
thus to sharpen the results of the previous work \cite{Me}. We also
note that the model in question could be easily generalized to
obtain the model of motion in $\R^d$.

\bigskip
Without loss of generality we assume that $X(0)=0$. Consider an
orthonormal basis of $\R^3$ such that for the acceleration $a$ it
is true that $a = (0,0,|a|)^{\top}$.

Our main result is the following
\begin{theo}
\label{Main} Suppose $0 < \alpha < 1$ and $|a| \neq 0$; then there
exist constants $c_1 > 0$ and $c_2, c_3 \ge 0$ such that for any
initial speed $v_0 \in \R^3$, in the space $\mathcal{C}
\bigl([0,1], \R^3 \bigr)$
$$Y_t(s):=\frac{X(st) - c_1 a st}{\sqrt{t}}
\stackrel{d}{\longrightarrow} Y(\cdot):=\left( \begin{array}{c} c_2 W_1(\cdot)\\ c_2 W_2(\cdot)\\
c_3 W_3(\cdot)\\ \end{array} \right), \qquad t \to \infty,$$ where
$W_1$, $W_2$, and $W_3$ are independent Wiener processes.
\end{theo}
\begin{rem*}
The constants $c_1, c_2$, and $c_3$ depend on the model's
parameters $a, \alpha$, and $\lambda$. The author failed to find
these constants in an explicit form.
\end{rem*}
\begin{rem*}
We can easily extend the model defined by the relations
\eqref{V_n+1 long} and \eqref{tau_n def} to get the model of
motion in random medium in $\R^d$. Indeed, let $\sigma_n$ be
uniformly distributed on $S^{d-1} \subset \R^d$, let $a \in \R^d$,
and let $F: \R^d \times \R_+ \to \R$ be defined as above. For this
model in $\R^d$, Theorem~\ref{Main} also holds; the proof is
almost the same as for $\R^3$. For the limit process, we have
$Y=(\,c_2 W_1, \, \dots \, , \, c_2 W_{d-1}, \, c_3 W_d)^{\top}$,
where $W_i$ are independent Wiener processes.

Moreover, for the motion in $\R^1$, Theorem~\ref{Main} is valid
for the similar model, where at a collision, the particle's speed
$v$ always changes to $-\alpha v$.
\end{rem*}

A nonrigorous explanation of our result could be found in
\cite{WE}.

Proving this theorem, we reduce the problem to some statements
about a certain Markov chain. The most difficult one is that for
this chain the functional central limit theorem (FCLT) holds. The
difficulties arise because the chain has uncountable and
noncompact state space; moreover, the chain does not satisfy the
Doeblin condition, thus the classical results (for instance, from
\cite{Doob}) are not applicable. We solve the problem using quite
recent results (see \cite{MT} and references therein), based on
stochastic analogues of Lyapunov's functions.

Therefore, our methods differ significantly from those of
\cite{R-Triolo}. Nevertheless, some similarities could be found on
a deep level. Indeed, the proof from \cite{R-Triolo} is based on
martingale theory, while the presented one follows from the FCLT
for Markov chains. Recall that the FCLT is proved by martingale
theory arguments.

\section{The deduction of the simplified Lorentz model} \label{Append}
In this section we consider the original Lorentz model, described in
Subsection~\ref{Motiv}. For this model, we find the distribution
of the time before the fist collision $\tau_0$ and the
distribution of the random normal $\nu_1$, which describes the
first collision. Then we use these results to derive the
simplified model, defined by \eqref{V_n+1 long} and \eqref{tau_n
def}.

Since the medium is assumed to be isotropic, it is natural to
suppose that the obstacles' centers form a Poisson point process,
with the control measure proportional to Lebesgue measure
$\lambda_3$. The only shortcoming of this assumption is that
obstacles can intersect. We denote this Poisson process by $\Pi$
and write its control measure in the form $(\pi \lambda r^2)^{-1}
\lambda_3$, where $r$ is the sum of radii of the particle and an
obstacle; $\lambda > 0$ is a parameter. We will see that
$\lambda$ signifies the mean free path of the particle.

In addition, we assume that for the initial speed $v_0$ the
following condition holds:
\begin{equation}
\label{in cond} |v_0^\bot|^2 > r|a|,
\end{equation}
where $v_0^\bot$ is the projection of $v_0$ on the orthogonal
complement of $a$.

\subsection{The distribution of $\tau_0$}
We start with the following notations: for a $v \in \R^3$, define
the unit hemisphere $$\mathcal{S}_v:=\bigl\{u \in \R^3: |u|=1, (u,
v) \ge 0 \bigr \};$$ for any $0 \le t_1 < t_2$, define the set
$$\mathcal{A}_{t_1,t_2}:=\left\{ v_0 s + \frac{a s^2}{2} +  r
\mathcal{S}_{v_0 + as}, \, s \in [t_1, t_2)\right\}.$$

Let us fix a $t > 0$. The inequality $\tau_0 > t$ is equivalent to
the absence of obstacles' centers in the set $\mathcal{A}_{0,t}$,
whence $$\P \bigl \{ \tau_0 > t \bigr \} = \P \bigl \{
\Pi(\mathcal{A}_{0,t})=0 \bigr \} = \exp \Bigl \{ - (\pi \lambda
r^2)^{-1} \lambda_3 ( \mathcal{A}_{0,t}) \Bigr \}.$$ We claim that
for volume of the set $\mathcal{A}_{0,t}$, which has the form of a
curved "cylinder",
\begin{equation}
\label{vol A_t} \lambda_3 ( \mathcal{A}_{0,t}) = \pi r^2 \int_0^ t
|v_0 + as|ds
\end{equation}
(the naive explanation is the following: the factor of the
integral is area of the "cylinder's" cross-section, and the
integral is the "cylinder's" length). Then, for the distribution
function of $\tau_0$,
\begin{equation}
\label{tau def} \P \bigl \{ \tau_0 > t \bigr \} =
e^{-\lambda^{-1}\int_0^ t |v_0 + as|ds}.
\end{equation}
Also, $\P \bigl \{ \tau_0 \le dt \bigr \} = \lambda^{-1}|v_0| dt +
o(dt)$, as $dt \to 0$, thus we see that the parameter $\lambda$
signifies the mean free path of the particle.

Let us prove \eqref{vol A_t}. Without loss of generality, assume
that $v_0^\bot=(0,|v_0^\bot|,0)^\top$; since $a=(0,0,|a|)^\top$,
the first coordinate of the trajectory of the curve $s \mapsto v_0
s + a s^2 /2$ is zero. Denote by $L(t):=\int_0^ t |v_0 + as|ds$
the length of this curve, and let $s \mapsto (0, \gamma_2(s),
\gamma_3(s))^\top$ be the natural parametrization of the curve
(that is a parametrization such that $\dot{\gamma_2}(s)^2 +
\dot{\gamma_3}(s)^2=1$ for any $s$). Then we can represent
$\mathcal{A}_{0,t}$ in the form
$$\mathcal{A}_{0,t}=\left\{ (0, \gamma_2(s), \gamma_3(s))^\top +
r \mathcal{S}_{(0, \dot{\gamma_2}(s), \dot{\gamma_3}(s))^\top}, \,
s \in [0, L(t))\right\},$$ and, finally, introducing the function
$$Q(s, \varphi, \theta):= \left(
\begin{array}{c}
0\\
\gamma_2(s)\\
\gamma_3(s)\\
\end{array} \right) + r \left(
\begin{array}{ccc}
1 & 0 & 0\\
0 & \dot{\gamma_2}(s) & -\dot{\gamma_3}(s)\\
0 & \dot{\gamma_3}(s) & \dot{\gamma_2}(s) \\
\end{array} \right) \left(
\begin{array}{c}
\sin \theta \cos \varphi\\
\cos \theta\\
\sin \theta \sin \varphi \\
\end{array} \right),$$ we have
$$\mathcal{A}_{0,t}= Q \bigl([0, L(t)),[0, 2\pi), [0, \pi/2] \bigr).$$

By simple, but tedious calculations, it follows from \eqref{in
cond} that for any $l >0$ the mapping $Q: [0, l) \times [0, 2\pi)
\times [0, \pi/2] \to \R^3$ is bijective. Further, using the
equalities $\dot{\gamma_2}(s)^2 + \dot{\gamma_3}(s)^2=1$ and
$\dot{\gamma_2}(s) \ddot{\gamma}_2(s) + \dot{\gamma_3}(s)
\ddot{\gamma}_3(s)=0$ (the second inequality is the derivative of
the first one), the reader can easily prove that for the Jacobian
of $Q$ it is true that $\mbox{Jac} \, Q(s, \varphi, \theta) =r^2
\sin \theta \cos \theta$. We finish the proof of \eqref{vol A_t}
integrating the Jacobian over the set $[0, L(t)) \times [0, 2\pi)
\times [0, \pi/2]$.

\medskip
There exists a very convenient representation of $\tau_0$. Let
$\eta$ be an exponential r.v. with mean $\lambda$; recall that $F:
\R^3 \times \R_+ \to \R$ is defined as the solution of the
equation
$$\int_0^{F(v, t)} |v + as|ds = t.$$ Since $F(v,
\cdot)$ is monotone and $F \left(v, \int_0^t |v + as|ds \right) =
t,$ we have
\begin{equation}
\label{tau = F}
\tau_0 \stackrel{d}{=} F(v_0, \eta).
\end{equation}

\subsection{The distribution of $\nu_1$}
Recall that at a collision the particle's speed $V_1=v_0+a\tau_0$
changes to $V_1 - (1+\alpha)(V_1 \cdot \nu_1) \nu_1$, where
$\nu_1$ is the inner unit normal to the first obstacle at the
point of collision. To simplify the notations, let us write $\nu$
instead of $\nu_1$. As in the previous subsection, we assume that
$v_0^\bot=(0,|v_0^\bot|,0)^\top$.

At first notice that $\nu \in \mathcal{S}_{V_1}$; this vector
could be defined by its spherical coordinates $(\varphi_{\nu},
\theta_{\nu}^{V_1})$, where the longitude $\varphi_{\nu} \in [0,
2\pi)$ is the angle between $(1,0,0)^\top$ and $\nu$, and the
latitude $\theta_{\nu}^{V_1} \in [0, \pi/2]$ is the angle between
$V_1$ and $\nu_1$. For any $0 \le t_1 < t_2$, $\varphi \in [0,
2\pi]$, and $\theta \in [0, \pi/2]$, define
$$\mathcal{A}_{t_1, t_2, \varphi, \theta}:=\left\{ v_0 s + \frac{a
s^2}{2} +  r \bigl\{ u \in \mathcal{S}_{v_0 + as}: \varphi_u
<\varphi, \theta_u^{v_0+as} <\theta \bigr\}, \, s \in [t_1,
t_2)\right\} \subset \mathcal{A}_{t_1,t_2}.$$

Then
\begin{eqnarray*}
\P \Bigl \{ \varphi_{\nu} <\varphi, \theta_{\nu} ^{V_1} < \theta
\, \Bigl. \Bigr | \, \tau_0 \in [t, t+dt) \Bigr \} &=& \frac{\P
\Bigl \{ \varphi_{\nu} <\varphi, \theta_{\nu} ^{V_1} < \theta ,
\tau_0 \in [t, t+dt) \Bigr \}}{ \P \bigl \{ \tau_0 \in [t, t+dt) \bigr \}} \\
&=& \frac{\P \bigl \{ \Pi(\mathcal{A}_ {0, t}) = 0, \Pi
(\mathcal{A}_{t, t + dt, \varphi, \theta }) = 1 \bigr\}} {\P \bigl
\{ \Pi(\mathcal{A}_ {0, t}) = 0, \Pi(\mathcal{A}_ {t, t+ dt}) = 1
\bigr\}} + o(1), \quad dt \to 0\\
\end{eqnarray*}
(the term $o(1)$ appears, because since $\lambda_3 (
\mathcal{A}_{t, t+dt} ) = O(dt)$, we have $\P \bigl \{
\Pi(\mathcal{A}_ {t, t+ dt}) \ge 2 \bigr\} = o(dt)$). The mapping
$Q$ from the previous subsection is bijective, therefore the sets
$\mathcal{A}_{0, t}$ and $\mathcal{A}_{t, t+dt}$ are disjoint.
Consequently, $$\P \Bigl \{ \varphi_{\nu} <\varphi, \theta_{\nu}
^{V_1} < \theta \, \Bigl. \Bigr | \, \tau_0 = t \Bigr \} = \lim
\limits_{dt \to 0} \frac{\P \bigl \{ \Pi (\mathcal{A}_{t, t + dt,
\varphi, \theta }) = 1 \bigr\}} {\P \bigl \{ \Pi(\mathcal{A}_ {t,
t+ dt}) = 1 \bigr\}} = \lim \limits_{dt \to 0} \frac{ \lambda_3
(\mathcal{A}_{t, t + dt, \varphi, \theta }) } { \lambda_3
(\mathcal{A}_ {t, t+ dt})};$$ but the numerator is the integral of
$\bigl | \mbox{Jac} \, Q \bigr|$ over $[L(t), L(t+dt)) \times [0,
\varphi) \times [0, \theta]$, the denominator is the integral of
$\bigl | \mbox{Jac} \, Q \bigr|$ over $[L(t), L(t+dt)) \times [0,
2 \pi) \times [0, \pi/2 ]$, and we get $$\P \Bigl \{ \varphi_{\nu}
<\varphi, \theta_{\nu} ^{V_1} < \theta \, \Bigl. \Bigr | \, \tau_0
= t \Bigr \} = \frac{\varphi \sin^2 \theta}{2\pi}.$$

We see that the distribution of $\nu_1$ is invariant under
rotations around $V_1$, and $$\P \bigl \{ \theta_{\nu_1}^{V_1} <
\theta \bigr \} = \sin^2 \theta;$$ moreover,
$\theta_{\nu_1}^{V_1}$ and $\tau_0$ are independent.

\medskip
Let us find a suitable representation of $\nu_1$. Suppose that
$\sigma$ is uniformly distributed on the unit sphere $S^2$ and is
independent of $V_1$ (and of $\tau_0$); let $\tilde{\nu}$ be the
unit vector directed along the bisectrix of the angle between
$\sigma$ and $V_1$. Then, for any fixed $V_1$, the conditional
distributions of $\nu_1$ and $\tilde{\nu}$ coincide. Indeed, both
of them are invariant under rotations around $V_1$, and $\P \bigl
\{ \theta_{\tilde{\nu}}^{V_1} < \theta \bigr \} = \P \bigl \{
\theta_{\sigma}^{V_1} < 2\theta \bigr \} = (1 - \cos 2 \theta)/2 =
\sin^2 \theta$. Finally, from $(V_1, \tilde{\nu}) \tilde{\nu} =
\frac{1}{2} (V_1 + |V_1| \sigma )$ it follows that $$V_1 -
(1+\alpha)(V_1 \cdot \nu_1) \nu_1 \stackrel{d}{=} V_1 - \frac{1 +
\alpha}{2} \bigl(V_1 + |V_1| \sigma \bigr).$$

\subsection{The simplified model}
We simplify the model by stating that we can apply the results of
two previous subsections to describe the particle's motion after
{\it each} collision (i.e., we can simply replace $v_0$ by the
speed after the collision). Thus we obtain the model defined by
\eqref{V_n+1 long} and \eqref{tau_n def}.

The main disparity between the models was discussed in
Subsection~\ref{my model}. In addition, the models differ, because
the distributions of $\tau_0$ and $\nu_1$ were derived under the
technical condition \eqref{in cond}; clearly, the variables $V_n -
\frac{1 + \alpha}{2} (V_n + |V_n| \sigma_n )$ and $v_0$ do not
have to satisfy it.

However, our simplifications look quite natural. Indeed, the
"disappearance" of obstacles after collisions in a sense means
that a particle never returns to already met obstacles. This is
reasonable if the obstacles are rare, because there is a drift in
the direction of $a$. It is also sensible to neglect \eqref{in
cond} if $r$ or $|a|$ is small; note that since distribution of
$V_n$ converges to $\pi_V$, the variables $\bigl (V_n - \frac{1 +
\alpha}{2} (V_n + |V_n| \sigma_n ) \bigr)^\bot $ converge to a
nondegenerate limit.

\section{Starting the proof of Theorem \ref{Main}}
\subsection{The position of the particle at the $n$th collision}
In this subsection we shall find a suitable representation for
$X_n := X(t_n)$, where $t_n := \sum_{i=0}^{n-1} \tau_i$ is the
moment of the $n$th collision; additionally, put $t_0:=0$.

Consider the new Markov chain
$$\Phi_n := \left( V_n \atop \sigma_n \right), \qquad n \in \N; \qquad
\Phi_0 := \left( v_0 \atop -v_0/|v_0| \right).$$ It will be
obvious later that the current initial condition describes the
particle with initial speed $v_0$ (at time zero happens the dummy
collision, which does not change the speed). Note that $V_n$ and
$\sigma_n$ are independent. Further, introducing the notations
$$\widehat{x} := v - \frac{1 + \alpha}{2} \bigl(v + |v| \sigma \bigr),
\qquad x = \left( v \atop \sigma \right) \in X := \R^3 \times
S^2,$$ we rewrite \eqref{V_n+1 long} and \eqref{tau_n def} as
$$V_{n+1} = \widehat{\Phi_n} + a F ( \widehat{\Phi_n}, \eta_n ).$$
Thus for the new chain it is true that $$\Phi_{n+1} = \left(
\widehat{\Phi_n} + a F ( \widehat{\Phi_n}, \eta_n ) \atop
\sigma_{n+1} \right), \qquad n \ge 0.$$

Let us agree to write coordinates of vectors of $\R^3$ using
superscripts; recall that $a^3 = |a|$. From the trivial equalities
$V^1_{n+1}=\widehat{\Phi_n}^1$, $V^2_{n+1}=\widehat{\Phi_n}^2$,
and $V^3_{n+1}=\widehat{\Phi_n}^3 + |a| \tau_n$ it follows that
$$ X_{n+1} = X_n + \left(
\begin{array}{c}
\widehat{\Phi_n}^1 \tau_n\\
\widehat{\Phi_n}^2 \tau_n\\
((V_{n+1}^3)^2 - (\widehat{\Phi_n}^3)^2 )/(2|a|)\\
\end{array}
\right)
= X_n + \frac{1}{|a| } \left(
\begin{array}{c}
V_{n+1}^1 V_{n+1}^3 - \widehat{\Phi_n}^1 \widehat{\Phi_n}^3\\
V_{n+1}^2 V_{n+1}^3 - \widehat{\Phi_n}^2 \widehat{\Phi_n}^3\\
((V_{n+1}^3)^2 - (\widehat{\Phi_n}^3)^2 )/2\\
\end{array}
\right),
$$
whence
\begin{equation}
\label{X_n=}
X_{n+1} = \frac{1}{|a| } \left( \begin{array}{c}
\widehat{\Phi_{n+1}}^1 \widehat{\Phi_{n+1}}^3 - \widehat{\Phi_0}^1 \widehat{\Phi_0}^3\\
\widehat{\Phi_{n+1}}^2 \widehat{\Phi_{n+1}}^3 - \widehat{\Phi_0}^2 \widehat{\Phi_0}^3\\
((\widehat{\Phi_{n+1}}^3)^2 - (\widehat{\Phi_0}^3)^2 )/2\\
\end{array}
\right) + \frac{1}{|a| } \sum_{i=1}^{n+1} \left( \begin{array}{c}
V_i^1 V_{i}^3 - \widehat{\Phi_i}^1 \widehat{\Phi_i}^3\\
V_i^2 V_{i}^3 - \widehat{\Phi_i}^2 \widehat{\Phi_i}^3\\
((V_i^3)^2 - (\widehat{\Phi_i}^3)^2 )/2\\
\end{array}
\right).
\end{equation}
Besides,
\begin{equation}
\label{t_n=} t_{n+1} = \sum_{i=0}^n \tau_i = \frac{1}{|a|}
\sum_{i=0}^n V_{i+1}^3 - \widehat{\Phi_i}^3 = \frac{1}{|a|} \bigl
(\widehat{\Phi_{n+1}}^3 - \widehat{\Phi_0}^3 \bigr) +
\frac{1}{|a|} \sum_{i=1}^{n+1} V_i^3 - \widehat{\Phi_i}^3.
\end{equation}
Finally, denoting $$f(x):=\frac{1}{|a| } \left( \begin{array}{c}
v^1 v^3 - \widehat{x}^1 \widehat{x}^3\\
v^2 v^3 - \widehat{x}^2 \widehat{x}^3\\
((v^3)^2 - (\widehat{x}^3)^2 )/2\\
\end{array}
\right), \quad  h(x):= \left( \begin{array}{c}
0\\
0\\
v^3 - \widehat{x}^3\\
\end{array}
\right), \qquad x = \left( v \atop \sigma \right) \in X,$$ from
\eqref{X_n=} and \eqref{t_n=} we get
\begin{equation}
\label{X_n-cat_n} X_n - c_1 a t_n = \frac{1}{|a| } \left(
\begin{array}{c}
\widehat{\Phi_n}^1 \widehat{\Phi_n}^3 - \widehat{\Phi_0}^1 \widehat{\Phi_0}^3\\
\widehat{\Phi_n}^2 \widehat{\Phi_n}^3 - \widehat{\Phi_0}^2 \widehat{\Phi_0}^3\\
((\widehat{\Phi_n}^3)^2 - (\widehat{\Phi_0}^3)^2 )/2\\
\end{array}
\right) - \frac{c_1}{|a|}a \bigl(\widehat{\Phi_n}^3 -
\widehat{\Phi_0}^3 \bigr) + \sum_{i=1}^n [f-c_1 h](\Phi_i).
\end{equation}
By definition, put $g:=f-c_1 h$ (the value of $c_1$ will be
defined later).

\subsection{The problem in terms of the Markov chain $\Phi_n$}
Between collisions the particle moves with constant acceleration,
thus the process $X(t)$ is defined by its values $X_n$ at the
points $t_n$. In fact, we can find the values of $X^1(t)$ and
$X^2(t)$ by linear interpolation and the values of $X^3(t)$ by
quadratic interpolation with the leading coefficient $|a|/2$.
Analogously, the process $Y_t(s)$ is defined by its values at the
points $t_n/t$, but for $Y_t^3(s)$ we shall use quadratic
interpolation with the leading coefficient $|a| t^{3/2}/2$.

Let $\widetilde{Y}_t(s)$ be the following process: at the points
$t_n/t$ put
$$\widetilde{Y}_t \Bigl(\frac{t_n}{t} \Bigr):= Y_t\Bigl(
\frac{t_n}{t} \Bigr) =\frac{X_n - c_1 a t_n}{\sqrt{t}}, \qquad n
\ge 0$$ and define the values at other points via linear
interpolation. The only difference between $Y_t(s)$ and
$\widetilde{Y}_t(s)$ is in the method of interpolation for the
third coordinate. It is easy to see that $Y^3_t(s) -
\widetilde{Y}^3_t(s) = |a| t^{3/2} (s-t_n) (s-t_{n+1})/2$ for $s
\in [t_n, t_{n+1}]$. Thus, denoting by $n(t)$ the (random) number
of collisions by the time $t$, for the norm $\| \cdot
\|_\mathcal{C}$ of the space $\mathcal{C} [0,1] = \mathcal{C}
\bigl([0,1], \R^3 \bigr)$ we have
\begin{equation}
\label{| |_1}
\| Y_t( \cdot) - \widetilde{Y}_t( \cdot) \|_\mathcal{C} \le \max
\limits_{0 \le k \le n(t)} \sup \limits_{t_k \le s \le t_{k+1}} |
Y_t(s) - \widetilde{Y}_t(s) | = \frac{|a|}{8 \sqrt{t}} \max
\limits_{0 \le k \le n(t)} \tau_k^2.
\end{equation}

Then, we introduce the process $Z_t(s)$, putting at the points
$t_n/t$
\begin{equation}
\label{Z_t def}
Z_t \Bigl(\frac{t_n}{t} \Bigr):=\frac{1}{\sqrt{t}} \sum_{i=1}^n
g(\Phi_i), \qquad n \ge 0
\end{equation}
and defining the values at other points via linear
interpolation. Trajectories of $Z_t(s)$ and $\widetilde{Y}_t(s)$
are piecewise linear and their points of interpolation have the
same $x$-coordinates (namely, $t_n/t$), therefore from
\eqref{X_n-cat_n} we have
\begin{equation}
\label{| |_2}
\| Z_t( \cdot) - \widetilde{Y}_t( \cdot) \|_\mathcal{C} \le
\frac{1}{\sqrt{t}} \max \limits_{1 \le k \le n(t)+1} \left\{
\frac{2.5}{|a|} \bigl(|\widehat{\Phi_k}|^2 + |\widehat{\Phi_0}|^2
\bigr) + c_1 \bigl(|\widehat{\Phi_k}| + |\widehat{\Phi_0}| \bigr)
\right\}.
\end{equation}

From \eqref{| |_1} and \eqref{| |_2} we see that for proving
Theorem~\ref{Main} it is sufficient to check that for all initial
conditions $\Phi_0=x \in X$, it is true that
\begin{equation}
\label{dif to 0}
\frac{1}{\sqrt{t}} \max \limits_{0 \le k \le n(t)} \tau_k^2
\stackrel{\P_x}{\longrightarrow} 0, \quad \frac{1}{\sqrt{t}} \max
\limits_{1 \le k \le n(t)+1} |\widehat{\Phi_k}|^2
\stackrel{\P_x}{\longrightarrow} 0, \qquad t \to \infty,
\end{equation}
and there exist constants  $c_1 >0 $ and $c_2, c_3 \ge 0$ such
that in the space $\mathcal{C}[0,1]$
\begin{equation}
\label{Z FCLT}
Z_t(\cdot) \stackrel{d}{\longrightarrow} Y(\cdot), \qquad t \to
\infty.
\end{equation}

Thus we must study the properties of the Markov chain $\Phi_n$ in
detail. The necessary facts from the Markov chain theory are stated in
Section~\ref{Svedeniya}. In Section~\ref{Research} we prove that
the chain $\Phi_n$ possess some useful properties. These
properties are used in Section~\ref{Proof} for proving \eqref{dif
to 0} and \eqref{Z FCLT}. The assertion \eqref{Z FCLT} is the most
difficult; we prove it applying the FCLT to the sequence
$g(\Phi_n)$.

The problem is reduced to studying properties of the Markov chain
$\Phi_n$.

\section{Basic facts on Markov chains} \label{Svedeniya}
The purpose of this section is to describe conditions under which
a Markov chain satisfies the law of large numbers (LLN) and
the FCLT. We also give a simple method for checking this
conditions. All the statements and definitions are taken from
\cite{MT}; in this section multiple references to this source are
omitted.

We begin with several notations. Consider a Markov chain $\Phi_n$,
with an arbitrary state space $X$ equipped with a locally compact,
separable, metrizable topology and Borel $\sigma$-field
$\mathcal{B}(X)$. Let $P(x, \cdot)$ be the transition function of
$\Phi_n$, let $P^n(x, \cdot)$ be the $n$-step transition function,
and let $\pi$ be the invariant measure of the chain (in the
considered situations, there exists a unique invariant probability
measure). Calculating expectations and probabilities, we indicate
the initial distribution of the chain, i.e., $\mathcal{L}(\Phi_0)$,
with subscripts. For example, $\P_x \bigl\{ \Phi_n \in A \bigr\}$
imply that $\mathcal{L}(\Phi_0) = \delta_x$ and $\E_\pi \Phi_n $
imply that $\mathcal{L}(\Phi_0) = \pi$. All the considered
functions are assumed to be measurable. Finally, by $P$ denote the
transition operator; recall that by definition $(Pf)(x)= \int_X
f(y) P(x,dy)$, for any functional $f: X \to \R$.

\subsection{Definitions}
A Markov chain is called {\it irreducible} if there exists a
nonzero measure $\mu$ on $\mathcal{B}(X)$ such that
$$\mu(A) > 0 \quad \Longrightarrow \quad \P_x \bigl\{ \exists
\, n \in \N : \Phi_n \in A \bigr\} > 0 , \qquad x \in X, A \in
\mathcal{B}(X);$$ any measure satisfying this condition is called
{\it irreducible measure} of the chain.

An irreducible chain is called {\it aperiodic} if there does not
exist a $d \ge 2$ and there do not exist disjoint sets $E_1,
\dots, E_{d} \in \mathcal{B}(X)$ such that
\\1) for all $x \in E_d$, $P(x, E_1)=1$, and for all $x \in E_i$, $P(x, E_{i+1})=1, \quad
i=1, \dots, d-1$;
\\2) $\mu \Bigl( X \, \Bigl. \Bigr \backslash \, \bigcup \limits_{i=1}^d E_i \Bigr) =0$ holds
for every irreducible measure $\mu$ of the chain \\(we modified
the definition from \cite{MT} using Proposition 4.2.2 and Theorem
5.4.4).

We say that a chain is {\it (weak) Feller } if the function $P(\,
\cdot \,, A)$ is lower semicontinuous for any open set $A \subset
X$.

Let $\mu$ be a signed measure on $\mathcal{B}(X)$, and let $f:X \to
[0, \infty)$ be a functional. We define the {\it $f$-norm} of
$\mu$ as
$$\| \mu \|_f:= \sup\limits_{g: |g| \le f} \int_X g d \mu = \int_X
f d |\mu|$$ (the inequality $|g| \le f$ is pointwise). The
$1$-norm is called the {\it total variation norm}; the notation
$\| \mu \|_1$ is replaced by $\| \mu \|$.

Let $P_1$ and $P_2$ be Markov transition functions, and let $U: X
\to [1, \infty)$ be a functional. By definition, put
$$|||P_1 - P_2|||_U := \sup\limits_{x \in X} \frac{\|P_1(x, \cdot)
- P_2(x, \cdot) \|_U } {U(x)}.$$

A Markov chain $\Phi_n$ is {\it ergodic} if there exists a measure
$\pi$ such that for any $x \in X$ it is true that $\|P^n(x, \cdot)
- \pi \| \to 0$, as $n \to \infty$; this yields that $\pi$ is a
unique invariant probability measure. A Markov chain is  {\it
U-uniformly ergodic} if there exists a measure $\pi$ such that
$|||P^n - \pi |||_U \to 0$, as $n \to \infty$ (we formally put
$\pi(x, \cdot):= \pi(\cdot)$). Note that if a chain is
$U$-uniformly ergodic, then it is $cU$-uniformly ergodic, for any
$c>1$. For irreducible aperiodic chains the $1$-uniform ergodicity
is equivalent to the well-known Doeblin condition (see Theorem
16.2.3).

Let $g: X \to \R$ be such that $g \in L^1 (\pi) = L^1(X,
\mathcal{B}(X), \pi)$. The functional equation (in unknown
$\bar{g}$)
\begin{equation}
\label{EqPuass}
\bar{g} - P \bar{g} = g - \int_{X} g d \pi
\end{equation}
is called the {\it Poisson equation}.

\subsection{Theorems}
\begin{theo}
\label{LLN}
Let $\Phi_n$ be an ergodic Markov chain, and let $g \in L^1( \pi)$.
Then for any initial condition $\Phi_0 = x \in X$
$$\lim\limits_{n \to \infty} \frac{1}{n} \sum_{i=1}^n g(\Phi_i)
= \int_{X} g d \pi , \qquad \P_x \mbox{{\it-a.s.}}$$
\end{theo}
\begin{proof}
Follows from Theorem 17.0.1.
\end{proof}
\begin{theo}\label{Erg criterium}
Let $\Phi_n$ be an irreducible, aperiodic, Feller Markov chain,
and let $\mbox{Int} \, (\mbox{supp} \, \mu) \neq \varnothing$ for
some irreducible measure $\mu$ of $\Phi_n$. Suppose that the
Foster-Lyapunov condition holds: there exist a functional $U: X
\to [1, \infty)$, a compact set $C \subset X$, and constants
$\beta, b > 0$ such that
\begin{equation}
\label{drift}
PU(x) - U(x) \le -\beta U(x) + b \I_C(x), \qquad x \in X.
\end{equation}
Then $\Phi_n$ is $U$-uniformly ergodic; moreover, $U \in L^1
(\pi)$.
\end{theo}
\begin{proof}
In view of Proposition 5.5.3 and Theorem 6.0.1, the first
statement follows from Theorems 15.0.1 and 16.0.1; the last one is
proved in Theorem 14.0.1.
\end{proof}
\begin{theo}
\label{FCLT}
Let $\Phi_n$ be a $U$-uniformly ergodic Markov chain, and let a
functional $g: X \to \R$ be such that $g^2 \le U$. Then there
exists a solution $\bar{g}$ of the Poisson equation
\eqref{EqPuass}; $\bar{g} \in L^2 (\pi)$; and the constant
\begin{equation}
\label{gamma}
\gamma_g^2 :=\int_X \left( \bar{g}^2 - (P \bar{g})^2 \right) d \pi
\ge 0
\end{equation}
is well defined. If $\int_{X} g d \pi=0$, then for any initial
condition $\Phi_0 = x \in X$, in the space $\mathcal{C}[0,1]$
\begin{equation}
\label{S_t FCLT}
S_t(s):=\frac{\sum_{i=1}^{[st]} g(\Phi_i) + (st - [st])
g(\Phi_{[st]+1}) }{\sqrt{t}} \stackrel{d}{\longrightarrow}
\sqrt{\gamma_g^2} W(\cdot), \qquad t \to \infty,
\end{equation}
where $W$ is a Wiener process.
\end{theo}
\begin{proof}
The existence of a solution of the Poisson equation easily follows
from Theorem 17.4.2. The well-posedness of the definition of
$\gamma_g^2$ (i.e., independence of the choice of the Poisson
equation's solution $\bar{g}$) follows from Proposition 17.4.1.
The Cauchy-Bunyakovskii-Schwarz inequality implies that
$\gamma_g^2 \ge 0$. Finally, the last statement is the combination
of Theorems 17.4.4 and 17.5.4. Although in \cite{MT} the processes
$S_t(s)$ are defined for positive integer $t$, in \eqref{S_t FCLT}
the convergence over $t \in \R$ simply follows from the
convergence over $t \in \N$.
\end{proof}

\section{Studying properties of the Markov chain $\Phi_n$} \label{Research}
We shall frequently use the following trivial inequalities: for
any $x =\left( v \atop \sigma \right) \in X$ it is true that
$$\alpha |v| \le |\widehat{x}| \le |v|;$$ recall that $\alpha \in (0, 1)$.

\subsection{Irreducibility, aperiodicity, and the Feller property}\label{psi_irr}
First, let us prove that for any  $\Phi_0 = x$ the values of $V_2
= \widehat{\Phi_1} + a F ( \widehat{\Phi_1}, \eta_1 )$ run over
the whole $\R^3$. It is sufficient to prove that
$\widehat{\Phi_1}$ runs through $\R^3 \setminus B(0,|\widehat{x}|)
= \bigl\{ u \in \R^3: |u| \ge |\widehat{x}| \bigr\}$, because for
any fixed $\widehat{\Phi_1}$ the values of $F ( \widehat{\Phi_1},
\eta_1 )$ run through $\R_+$. We certainly assume that $\eta_i$
and $\sigma_i$ run over the whole $\R_+$ and $S^2$ respectively.

At a collision, the speed $v \in \R^3$ changes to $\tilde{v}:= v -
(1+\alpha)(v, \nu) \nu$ (we temporary use the old representation).
The reader will easily check that the inverse transformation is $v
= \tilde{v} - (1+\alpha^{-1}) (\tilde{v}, \nu) \nu$, where $\nu$
is the same as in the direct transformation. The values of $\nu$
run through $\mathcal{S}_v=\bigl\{ u \in \R^3: |u|=1, \, (v, u)
\ge 0 \bigr\}$; consequently, we have $(\tilde{v}, \nu) \le 0$,
that is, $\nu \in \mathcal{S}_{-\tilde{v}}$. Thus the speed {\it
after} a collision could be equal to a $v \in \R^3$ iff the speed
{\it before} this collision is contained in the set $v^{-1}:=
\bigl\{ u \in \R^3: u = v - (1+\alpha^{-1})(v, \nu) \nu, \, \nu
\in \mathcal{S}_{-v} \bigr\}$. As above, define $\sigma$ as a unit
vector such that $\nu$ is directed along the bisectrix of the
angle between $\sigma$ and $-v$. Then $\sigma$ runs over $S^2$,
whence $v^{-1}= \bigl\{ u \in \R^3: u = v +
\frac{1+\alpha^{-1}}{2} (|v|\sigma  - v ), \, \sigma \in S^2
\bigr\}$.

We shall prove that $\widehat{\Phi_1}$ runs over $\R^3 \setminus
B(0,|\widehat{x}|)$; recall that $\Phi_1 = \left( \widehat{x} + a
F(\widehat{x}, \eta_0) \atop \sigma_1 \right)$. Let us show that
for any $v \in \R^3 \setminus B(0,|\widehat{x}|)$ the set $v^{-1}
\bigcap \bigl\{ \widehat{x} + as, \, s \ge 0 \bigr\}$ is nonempty.
This is equivalent to the existence of a $s \ge 0$ such that
$\bigl|\widehat{x} + as - \frac{1-\alpha^{-1}}{2}v \bigr| =
\frac{1+\alpha^{-1}}{2}|v|$. But the left-hand side continuously
depends on $s$ and increases for large $s$. Since $|v| \ge
|\widehat{x}|$, the value of the left-hand side at $s=0$ is not
greater than $\frac{1+\alpha^{-1}}{2}|v|$. Therefore the required
$s$ exists.

Moreover, the distribution of $\widehat{\Phi_1}$ has a density,
and this density is positive on $\R^3 \setminus
B(0,|\widehat{x}|)$. This fact is intuitively clear, because
$F(\widehat{x}, \eta_0)$ has a positive density on $\R_+$ (see
\eqref{tau def} and \eqref{tau = F}) and $\sigma_1$ has a positive
density on $S^2$. The formal proof, whose main part is to
calculate the Jacobian of the appropriate transformation, is
omitted.

Now it is obvious that the distribution of $V_2 = \widehat{\Phi_1}
+ a F ( \widehat{\Phi_1}, \eta_1 )$ has a positive density on
$\R^3$. That is why the chain $\Phi_n$ is irreducible and
$\lambda_3 \otimes U_{S^2}$ is an irreducible measure; here
$\lambda_3$ is the Lebesgue measure on $\R^3$ and $U_{S^2}$ is the
uniform distribution on $S^2$. Indeed, $\Phi_2 = \left( V_2 \atop
\sigma_2 \right)$, and the variables $V_2$ and $\sigma_2$ are
independent.

\medskip
We just proved that for any initial condition $\Phi_0 = x$ the
distribution of $\Phi_2$, i.e., $P^2(x, \cdot)$, has positive
density with respect to $\lambda_3 \otimes U_{S^2}$. Hence for
{\it any} $x \in X$ the measures $P^2(x, \cdot)$ and $\lambda_3
\otimes U_{S^2}$ are equivalent. Thus, by simple arguments, the
chain $\Phi_n$ is aperiodic.

\medskip
To prove that $\Phi_n$ is a Feller chain, it is sufficient to
check that for any open set $A \subset \R^3$ the function $P(\,
\cdot \,,A \times S^2)$ is lower semicontinuous. Indeed, from
\eqref{tau def} it follows that
$$P(x, A \times S^2) = \int_0^\infty \I_A(\widehat{x}+at)
\lambda^{-1} |\widehat{x} + at| e^{-\lambda^{-1}\int_0^ t
|\widehat{x} + as|ds} dt,$$ and applying the Fatou lemma, we
obtain lower semicontinuity.

\subsection{$U$-uniform ergodicity}
\begin{lem}
\label{exp moment} Let $\eta$ be an exponential r.v. with mean
$\lambda$; then for any $c \in \R$ $$\sup \limits_{v \in \R^3} \E
e^{c F(v, \eta)} < \infty.$$
\end{lem}

\begin{cor*}
For any initial condition $\Phi_0=x \in X$, the variables $\tau_n$
have exponential moments of any order. Moreover,
\begin{equation}
\label{tau exp b}
\sup \limits_n \E_x e^{c \tau_n} < \infty.
\end{equation}
\end{cor*}
\begin{proof}[{\bf Proof of Corollary}]
Since $\tau_n = F(\widehat{\Phi_n}, \eta_n)$ and
$\widehat{\Phi_n}$ is independent of $\eta_n$, the proof is
obvious.
\end{proof}
\begin{proof}[{\bf Proof of Lemma~\ref{exp moment}}]
We only consider the nontrivial case $c > 0$. Take an $s > 0$ and
a $v \in \R^3$. Since $\bigl| -a|v|/|a| + as \bigr| = \bigl|
|a|s-|v| \bigr| \le |v + as|$, we see that for any $t>0$
$$\int_0^t \bigl|-a|v|/|a| + as \bigr| ds \le \int_0^t |v + as|
ds. $$ Then, using the definition of $F$, we have $F(v, \cdot)
\le F(-a|v|/|a|, \cdot)$ and therefore $\P \bigl\{ F(v, \eta) > t
\bigr\} \le \P \bigl\{ F(-a|v|/|a|, \eta) > t \bigr\}$. The right-hand
side could be easily calculated, and
\begin{eqnarray*}
\E e^{c F(v, \eta)} &=& -\int_0^{\infty} e^{ct} d \, \P \{ F(v,
\eta)
> t \} \\
&=& c \int_0^{\infty} e^{ct} \P \{ F(v, \eta) > t \} dt - 1 \\
&\le& c \int_0^{\infty} e^{ct} \P \{ F(-a|v|/|a|, \eta) > t \} dt
- 1\\
&=& c \int_0^{\infty} e^{ct} e^{- \lambda^{-1} \int_0^t \left|
|a|s-|v| \right| ds } dt - 1 \\
&=& c \int_0^{|v|/|a|} e^{ct - \lambda^{-1}(|v|t-|a|t^2/2)} dt + c
\int_{|v|/|a|}^{\infty} e^{ct - \lambda^{-1}
(|a|t^2/2 - |v|t + |v|^2/|a|)} dt - 1.\\
\end{eqnarray*}

We now estimate the integrals. For the first one, use the
following inequality: if $|v| \ge 4 \lambda c$, then $ct -
\lambda^{-1} (|v|t-|a|t^2/2) \le -ct$, for all $t \in [0,|v|/|a|]$.
Therefore if $|v| \ge 4 \lambda c$, then
$$\int_0^{|v|/|a|} e^{ct - \lambda^{-1}(|v|t-|a|t^2/2)} dt \le
\int_0^{|v|/|a|} e^{-ct} dt < \int_0^{\infty} e^{-ct} dt <\infty$$
and this bound does not depend on $v$. If $|v| <  4 \lambda c$,
then $$\int_0^{|v|/|a|} e^{ct - \lambda^{-1}(|v|t-|a|t^2/2)} dt <
\int_0^{4\lambda c/|a|} e^{ct +\lambda^{-1} |a|t^2/2} dt
<\infty.$$

Estimating the second integral, we apply the following: if $|v|
\ge 5 \lambda c$, then $ct - \lambda^{-1} (|a|t^2/2 - |v|t +
|v|^2/|a|) \le -ct$. Thus if  $|v| \ge 5\lambda c$, then
$$\int_{|v|/|a|}^{\infty} e^{ct - \lambda^{-1}
(|a|t^2/2 - |v|t + |v|^2/|a|)} dt \le \int_{|v|/|a|}^{\infty}
e^{-ct} dt < \int_0^{\infty} e^{-ct} dt < \infty.$$ If $|v| <
5\lambda c$, then we have $$\int_{|v|/|a|}^{\infty} e^{ct -
\lambda^{-1} (|a|t^2/2 - |v|t + |v|^2/|a|)} dt < \int_0^{\infty}
e^{6ct - \lambda^{-1} |a|t^2/2 } dt < \infty.$$
\end{proof}

To prove the $U$-uniform ergodicity of $\Phi_n$, we apply
Theorem~\ref{Erg criterium}. Take a $c > 0$ and check that the
Foster-Lyapunov condition \eqref{drift} holds for $U(x) := e^{c
|\widehat{x}|}$.

For the transition operator,
$$PU(x) = \E_x U(\Phi_1) = \E_x e^{c | \widehat{\Phi_1}|} =
\E_x e^{c | V_1 - \frac{1 + \alpha}{2} (V_1 + |V_1| \sigma_1 )|
}.$$ We define the function $$\gamma(|v|):=\int_{S^2} e^{c \left|
v - \frac{1 + \alpha}{2} (v + |v|\zeta) \right| - c|v|} d
U_{S^2}(\zeta), \qquad v \in \R^3,$$ which is obviously monotone
and $\gamma (|v|) \to 0$ as $|v| \to \infty.$ Because $\sigma_1$
is independent of $V_1$ and of $\Phi_0$,
\begin{eqnarray*}
PU(x) &=& \E_x \gamma (|V_1|) e^{c |V_1|} \\
&\le& \E_x \gamma (|V_1|) e^{c | \widehat{\Phi_0}| } e^{c
|a| \tau_0} \\
&\le& \E_x \left(\gamma (|\widehat{\Phi_0}/2|) \I_{ \{ |a| \tau_0
\le |\widehat{\Phi_0}|/2 \}}  + \gamma(0) \I_{ \{ |a| \tau_0 >
|\widehat{\Phi_0}|/2 \} } \right) e^{c |
\widehat{\Phi_0}| } e^{c |a| \tau_0} \\
&\le& \left( \gamma (| \widehat{x}/2|) \E e^{c |a| F(\widehat{x},
\eta_0)} + \E e^{c |a| F(\widehat{x}, \eta_0)} \I_{ \{ |a|
F(\widehat{x}, \eta_0) > |\widehat{x}|/2 \} }\right) e^{c
|\widehat{x}|}\\
&\le& \left(\gamma (\alpha |v|/2) \E e^{c |a| F(\widehat{x},
\eta_0)} + \E e^{c |a| F(\widehat{x}, \eta_0)} \I_{ \{ |a|
F(\widehat{x}, \eta_0) > \alpha |v|/2 \} }\right) U(x),
\end{eqnarray*}
where as usual $x = \left( v \atop \sigma \right)$. It follows
from Lemma~\ref{exp moment} that the factor of $U(x)$ tends to
zero as $|v| \to \infty$, whence for any $\beta \in (0,1)$ there
exists an $R>0$ such that
$$PU(x) - U(x) \le - \beta U(x), \qquad x \notin C_R:=B(0,R)
\times S^2.$$ Clearly, for some $b > 0$ $$PU(x) - U(x) \le - \beta
U(x) + b \I_{C_R}(x), \qquad x \in X.$$ Thus the condition
\eqref{drift} holds and, consequently, for any $c > 0$ the Markov
chain $\Phi_n$ is $e^{c |\widehat{x}|}$-uniformly ergodic.

\subsection{The invariant measure} \label{inv mes}
By definition of $U$-uniform ergodicity, there exists a unique
invariant measure $\pi$ of the chain. Since for every $n$ the
measure $P^n(x, \cdot)$ is a product of some measure on
$\mathcal{B}(\R^3)$ and $U_{S^2}$, for the limit we also have
$$\pi = \pi_V \otimes U_{S^2},$$ where $\pi_V$ is a probability
measure on $\mathcal{B}(\R^3)$. By the reasons of symmetry,
$\pi_V$ is invariant under rotations around the third coordinate
axis.

Further, we claim that the measure $\pi_V$ has a density. Indeed,
in Subsection~\ref{psi_irr} we proved that $P^2(x, \cdot)$ has a
density; moreover, it could be shown that for all $n \ge 2$ the
measures $P^n(x, \cdot)$ have densities, that is, $P^n(x, \cdot)
\prec \lambda_3 \otimes U_{S^2}$. Thus, passing to the limit,
$\pi=\pi_V \otimes U_{S^2} \prec \lambda_3 \otimes U_{S^2}$ and
$\pi_V \prec \lambda_3$.

Theorem~\ref{Erg criterium} implies that $e^{c |\widehat{x}|} \in
L^1(\pi)$ for any $c$, and we can prove the following
\begin{pred}
For any initial condition $\Phi_0=x \in X$, the variables $V_n$
have exponential moments of any order. Moreover,
\begin{equation}
\label{V exp b}
\sup \limits_n \E_x e^{c |V_n|} < \infty.
\end{equation}
\end{pred}
\begin{proof}
We begin with $|V_n|= |\widehat{\Phi_{n-1}} + a \tau_{n-1} | \le
|V_{n-1} | + |a| \tau_{n-1} \le \dots \le |\widehat{x}| +|a|
\sum_{i=0}^{n-1} \tau_i$. Then, via the H\"{o}lder inequality and
Lemma~\ref{exp moment},
$$\E_x e^{c |V_n|} \le e^{c |\widehat{x}|} \E_x \prod_{i=0}^{n-1}
e^{c |a| \tau_i} \le e^{c |\widehat{x}|} \prod_{i=0}^{n-1} (\E_x
e^{c |a| n \tau_i}) ^{1/n} \le e^{c |\widehat{x}|} \sup \limits_{v
\in \R^3} \E e^{c |a| n F(v, \eta)}< \infty,$$ hence the
exponential moments exist. To prove \eqref{V exp b}, we combine
the trivial inequality $e^{c|u|} \le e^{c \alpha^{-1}
|\widehat{y}| }$, for $y = \left( u \atop \rho \right) \in X$, and
the definition of $e^{c \alpha^{-1} |\widehat{y}|}$-uniform
ergodicity of $\Phi_n$:
$$\lim \limits_{n \to \infty} \E_x e^{c |V_n|} = \lim \limits_{n
\to \infty} \int_X e^{c|u|} P^n(x, dy) = \int_X e^{c|u|} d \pi(y)
\le \int_X e^{c \alpha^{-1} |\widehat{y}|} d \pi(y) < \infty.$$
\end{proof}

\section{Finishing the proof of Theorem~\ref{Main}} \label{Proof}
Recall that we must prove \eqref{dif to 0} and \eqref{Z FCLT}. Let
us start with the following lemmas.
\begin{lem}
\label{a.s. bounds} For any initial condition $\Phi_0=x \in X$,
$$\qquad \tau_n = O(\log n), \qquad | V_n| = O(\log n) , \qquad \P_x \mbox{{\it
-a.s.}}$$
\end{lem}
\begin{proof}
Applying the Chebyshev inequality and then using \eqref{tau exp
b}, we have $$\P_x\bigl\{ \tau_n > 2 \log n \bigr\} \le \frac{\E_x
e^{\tau_n}} {n^2} \le \frac{1}{n^2} \sup \limits_n \E_x e^{c
\tau_n} .$$ Thus the first statement immediately follows from the
Borel-Cantelli lemma. Similarly, we prove the second statement via
\eqref{V exp b}.
\end{proof}

\begin{lem}
\label{t_n/n}
There exists a $c_4>0$ such that for any initial condition
$\Phi_0=x \in X$,
$$\lim \limits_{n \to \infty} \frac{t_n}{n} = c_4, \qquad \P_x
\mbox{{\it -a.s.}}$$
\end{lem}
\begin{proof}
Recalling \eqref{t_n=} and the introduced notations, we see that
it is sufficient to prove
$$\lim \limits_{n \to \infty} \frac{\widehat{\Phi_n}^3}{n} = 0,
\qquad \P_x \mbox{{\it -a.s.}};$$ the existence of a $c_4$ such
that
$$\lim \limits_{n \to \infty} \frac{1}{n |a|}\sum_{i=1}^n h^3(\Phi_i)
= c_4, \qquad \P_x \mbox{{\it -a.s.}};$$ and positiveness of
$c_4$.

Since $|\widehat{\Phi_n}| \le |V_n|$, from Lemma~\ref{a.s. bounds}
we immediately obtain the first statement. Further, we can apply
Theorem~\ref{LLN} to prove the second statement, because $|h^3(x)|
= |v^3 - \widehat{x}^3 |\le (1 + \alpha^{-1})|\widehat{x}| < (1 +
\alpha^{-1}) e^{|\widehat{x}|} \in L^1(\pi)$ and thus $h^3 \in
L^1(\pi)$. By definition, put $c_4:= |a|^{-1} \int_X h^3 d \pi$.
The proof of positiveness of $c_4$, which is quite simple, could
be found in \cite{Me}.
\end{proof}

\subsection{Proof of \eqref{dif to 0}}
For the r.v. $n(t)$, which denotes the number of collisions by the
time $t$, it is true that
\begin{equation}
\label{n(t)}
\lim\limits_{t \to \infty} n(t) = \infty, \qquad \P_x \mbox{{\it
-a.s.}}
\end{equation}
To prove this, assume the converse. Then we can find a $k \ge 0$
such that with nonzero probability the particle collides with
obstacles only $k$ times. Thus the probability of $\tau_k =
\infty$ is nonzero that contradicts with the existence of
exponential moments.

By Lemma~\ref{t_n/n} and \eqref{n(t)}, $$\lim \limits_{t \to
\infty} \frac{t_{n(t)}}{n(t)} = c_4, \qquad \P_x \mbox{{\it
-a.s.}},$$ but since $t_{n(t)} \le t < t_{n(t)+1}= t_{n(t)} +
\tau_{n(t)}$, from Lemma~\ref{a.s. bounds} we have
\begin{equation}
\label{t/n(t)}
\lim \limits_{t \to \infty} \frac{t}{n(t)} = c_4, \qquad \P_x
\mbox{{\it -a.s.}}
\end{equation}
Applying Lemma~\ref{a.s. bounds} once again (recall that
$|\widehat{\Phi_k}| \le |V_k|$), we prove \eqref{dif to 0} (to be
precise, we prove $\P_x$-a.s. convergence, which is much
stronger).

\subsection{Definition of the constants $c_1, c_2, c_3$}
We put $$c_1:= c_4^{-1} |a|^{-1} \int_X f^3 d \pi;$$ for the quite
tedious proof of positiveness of $c_1$, see \cite{Me}.

Now we can easily check that for $g=f-c_1 h$ it is true that
$\int_X g d \pi=0$. In fact, for the first and the second
coordinates, this follows from simple calculations, where the
representation $\pi = \pi_V \otimes U_{S^2}$ and symmetry of
$\pi_V$ are used. For the third coordinate, we apply the equality
$c_4 = |a|^{-1}\int_X h^3 d \pi $.

Let us define $c_2$ and $c_3$. By $|g(x)| \le 2.5 |a|^{-1} (1 +
\alpha^{-1}) |\widehat{x}|^2 + c_1 (1 +
\alpha^{-1})|\widehat{x}|$, there exists a $c>1$ such that
$|g(x)|^2 \le c e^{|\widehat{x}|}$. Thus (the chain $\Phi_n$ is $c
e^{|\widehat{x}|}$-uniformly ergodic) the functionals  $g^1$,
$g^2$, and  $g^3$ satisfy the conditions of Theorem~\ref{FCLT},
and there exist solutions $\overline{g^1}, \overline{g^2},
\overline{g^3} \in L^2(\pi)$ of the Poisson equations. Define
$$c_2:=\sqrt{c_4^{-1} \gamma^2_{g_2}}, \qquad c_3:=\sqrt{c_4^{-1}
\gamma^2_{g_3}},$$ and also $\bar{g}:= (\overline{g^1},
\overline{g^2}, \overline{g^3})^\top$,
$$K:=\int_X \left( \bar{g} \, \bar{g}^{\top} - (P
\bar{g})(P\bar{g})^{\top} \right) d \pi.$$ In \cite{Me} we proved
(using the axial symmetry of $\pi_V$) that the matrix $K$ is
diagonal and $\gamma^2_{g_1}= \gamma^2_{g_2}$, thus
$$K = c_4 \left(
\begin{array}{ccc}
c_2^2 & 0 & 0\\
0 & c_2^2 & 0\\
0 & 0 & c_3^2\\
\end{array} \right).$$

\subsection{Proof of \eqref{Z FCLT}}
In this subsection the following proposition plays the key role.
\begin{pred}
For the processes $S_t(s)$, defined in \eqref{S_t FCLT},
$$S_t(\cdot) \stackrel{d}{\longrightarrow} \sqrt{c_4} Y(\cdot),
\qquad t \to \infty.$$
\end{pred}
\begin{proof}
It is sufficient to show that for any $u \in \R^3$
$$(S_t(\cdot),u) \stackrel{d}{\longrightarrow} \sqrt{c_4}
(Y(\cdot),u), \qquad t \to \infty.$$ On the one hand, $$\sqrt{c_4}
(Y(\cdot),u) \stackrel{d}{=} \sqrt{c_4 \bigl((c_2 u^1)^2 + (c_2
u^2)^2 + (c_3 u^3)^2 \bigr)} W(\cdot) = \sqrt{(Ku, u)} W(\cdot).$$
On the other hand, the functional $(g, u)$ satisfies the
conditions of Theorem~\ref{FCLT}, thus
$$(S_t(s), u)=\frac{\sum_{i=1}^{[st]} (g,u)(\Phi_i) + (st - [st])
(g,u)(\Phi_{[st]+1}) }{\sqrt{t}} \stackrel{d} {\longrightarrow}
\sqrt{\gamma_{(g,u)}^2} W(\cdot), \qquad t \to \infty.$$ And
since, obviously, $\overline{(g, u)} =(\bar{g}, u)$, from
\eqref{gamma} we get $$\gamma^2_{(g, u)} = \int_{X} \left(
(\bar{g}, u)^2 - (P (\bar{g}, u))^2 \right) d \pi = \int_{X}
\left( (\bar{g}, u)^2 - (P \bar{g}, u)^2 \right) d \pi = (Ku,
u).$$
\end{proof}

Let us put $\widetilde{S}_t(s):= S_t(c_4^{-1}s)$; then
$S_t(c_4^{-1}s)=\sqrt{c_4^{-1}} S_{c_4^{-1}t}(s)$, and in the
space $\mathcal{C}[0,1]$ (and, moreover, in $\mathcal{C}[0,l]$,
for every $l>0$)
\begin{equation}
\label{S tilde ->}
\widetilde{S}_t(\cdot) \stackrel{d}{\longrightarrow} Y(\cdot),
\qquad t \to \infty.
\end{equation}
Therefore we will prove \eqref{Z FCLT} if we show that
\begin{equation}
\label{| |_3}
\|Z_t(\cdot )- \widetilde{S}_t(\cdot) \|_\mathcal{C}
\stackrel{\P_x} {\longrightarrow} 0, \qquad t \to \infty.
\end{equation}

For this purpose, introduce the process $u_t(s)$, putting at the
points $t_n/t$ $$\quad u_t \Bigl(\frac{t_n}{t} \Bigr):= \frac{c_4
n}{t}, \qquad n \ge 0$$ and defining the values at other points
via linear interpolation. Using the definitions of $Z_t(s)$ and
$S_t(s)$, i.e., \eqref{Z_t def} and \eqref{S_t FCLT}, we have $Z_t
\bigl( t_n/t \bigr) = S_t \bigl( n/t \bigr)$. But $S_t \bigl( n/t
\bigr) = \widetilde{S}_t \bigl(c_4n/t \bigr) = \widetilde{S}_t
\bigl( u_t \bigl(t_n/t \bigr) \bigr)$, whence at the points
$t_n/t$ the equality $Z_t \bigl(t_n/t \bigr) = \widetilde{S}_t
\bigl( u_t \bigl(t_n/t \bigr) \bigr)$ holds. However,
$$Z_t (s) = \widetilde{S}_t (u_t(s))$$ is true for {\it every} $s$ !
In fact, trajectories of $Z_t(s)$ and $\widetilde{S}_t (u_t(s))$
are piecewise linear (for the last one, as a composition of
piecewise linear functions) and their points of interpolation have
the same $x$-coordinates (namely, $t_n/t$). As we saw before, the
values at these points coincide.

We see that $Z_t(s)$ is obtained from $\widetilde{S}_t (s)$ by the
random change of time. Suppose $\| u_t(\cdot)\|_\mathcal{C} \le
2$; then
$$\|Z_t(\cdot )- \widetilde{S}_t(\cdot) \|_\mathcal{C}  =
\|\widetilde{S}_t( u_t(\cdot))- \widetilde{S}_t(\cdot)
\|_\mathcal{C} \le \omega_{\widetilde{S}_t \bigl. \bigr |_{[0,2]}}
\bigl(\| u_t(\cdot) - \mbox{id}\|_\mathcal{C} \bigr),$$ where
$\omega$ is the modulus of continuity, $\widetilde{S}_t \bigl.
\bigr |_{[0,2]}$ is the restriction of $S_t(s)$ to $[0,2]$. Hence
for any $1 > \delta >0$ and $\varepsilon>0$
\begin{eqnarray*}
\P_x \bigl\{ \|Z_t(\cdot )- \widetilde{S}_t(\cdot) \|_\mathcal{C}
\ge \varepsilon \bigr\} &\le& \P_x \bigl\{ \| u_t(\cdot) -
\mbox{id}\|_\mathcal{C} \ge \delta \bigr\} + \P_x \bigl\{
\omega_{\widetilde{S}_t \bigl. \bigr |_{[0,2]}} (\delta) \ge
\varepsilon \bigr\}\\
&\le& \P_x \bigl\{ \| u_t(\cdot) - \mbox{id}\|_\mathcal{C} \ge
\delta \bigr\} + \sup \limits_{t > 0} \P_x \bigl\{
\omega_{\widetilde{S}_t \bigl. \bigr |_{[0,2]}} (\delta) \ge
\varepsilon \bigr\}.\\
\end{eqnarray*}
Let us proceed to the limit as $t \to \infty$ and then proceed to
the limit as $\delta \to 0$. Now it is obvious that \eqref{| |_3}
holds if for any $1 > \delta > 0$
\begin{equation}
\label{|id - u_t|}
\lim \limits_{t \to \infty} \P_x \bigl\{ \| u_t(\cdot) -
\mbox{id}\|_\mathcal{C} \ge \delta \bigr\} = 0
\end{equation}
and for any $\varepsilon > 0$
\begin{equation}
\label{sup omega}
\lim \limits_{\delta \to 0} \,\sup \limits_{t >
0} \P_x \bigl\{ \omega_{\widetilde{S}_t \bigl. \bigr |_{[0,2]}}
(\delta) \ge \varepsilon \bigr\} = 0.
\end{equation}

At first we prove \eqref{|id - u_t|}. Writing $u_t(s)$ in the
explicit form, we have $$\|u_t(\cdot) - \mbox{id} \|_\mathcal{C} =
\sup \limits_{0 \le s \le 1} \Bigl| \frac{c_4 n(st)}{t} +
\frac{st-t_{n(st)}}{\tau_{n(st)}} \cdot \frac{c_4}{t} - s \Bigr|
\le \sup \limits_{0 \le s \le 1} \Bigl | \frac{c_4 n(st) - st}{t}
\Bigr | + \frac{c_4}{t};$$ thus, by \eqref{t/n(t)},
$$\lim \limits_{t \to \infty}\|u_t(\cdot) - \mbox{id}
\|_\mathcal{C} = 0, \qquad \P_x \mbox{{\it -a.s.}},$$ which is
much stronger than \eqref{|id - u_t|}.

It remains to check \eqref{sup omega} to complete the proof of
Theorem~\ref{Main}. The family of probability measures $\{ \P_x
\circ \widetilde{S}_t \bigl. \bigr |_{[0,2]}^{-1} \}_{t>0}$ on
$\mathcal{B} \bigl (\mathcal{C}[0, 2] \bigr)$ is relatively weakly
compact. This follows from \eqref{S tilde ->} and from $\P_x$-a.s.
continuity of $\| \widetilde{S}_t \|_{\mathcal{C} [0,2]}$ in $t
\in [0, \infty)$. The space $\mathcal{C}[0,2]$ is a Polish space,
thus the relatively weakly compact family $\{ \P_x \circ
\widetilde{S}_t \bigl. \bigr |_{[0,2]}^{-1} \}_{t>0}$ is tight. By
the well-known fact (see \cite{Bil}) about tight families of
probability measures on $\mathcal{B} \bigl (\mathcal{C}[0, 2]
\bigr)$, for any $\varepsilon >0$ $$\lim \limits_{\delta \to 0}
\,\sup \limits_{t > 0} \P_x \circ \widetilde{S}_t \bigl. \bigr
|_{[0,2]}^{-1} \Bigl\{p: \, \omega_p (\delta) \ge \varepsilon
\Bigr\} = \lim \limits_{\delta \to 0} \,\sup \limits_{t > 0} \P_x
\Bigl\{\omega_{\widetilde{S}_t \bigl. \bigr |_{[0,2]}} (\delta)
\ge \varepsilon \Bigr\} = 0.$$

\section*{Acknowledgements}
The author is deeply grateful to his advisor M.A. Lifshits, whose
remarks and numerous advices significantly improved this paper.
The author also acknowledges the anonymous referee for helpful
comments.

\bigskip
\begin{tabular}{>{\footnotesize} l}
Vladislav V. Vysotsky \\
Department of Probability Theory and Mathematical Statistics\\
Faculty of Mathematics and Mechanics\\
St.-Petersburg State University\\
Bibliotechnaya pl., 2 \\
Stariy Peterhof, 198504\\
Russia\\
E-mail: vysotsky@vv9034.spb.edu\\
\end{tabular}

\end{document}